\documentstyle[12pt]{article}
\def\ba{\begin{eqnarray}}
\def\ea{\end{eqnarray}}
\def\lb{\label}
\def\be{\begin{equation}}
\def\ee{\end{equation}}

\def\R{\hat{R}}
\def\F{\hat{F}}
\def\id{{\mit I}}

\newcommand{\tr}{\mathop{\mathrm{Tr}}\nolimits}
\newcommand{\Rtr}{\mathop{\mathrm{Tr}_{_{\!\R}}}\nolimits}

\begin{document}

\vspace{1cm}

\begin{center}
{\Large \bf Cayley-Hamilton-Newton identities and}

\vskip .6cm

{\Large \bf quasitriangular Hopf algebras}
\vskip 1cm
{A.P. Isaev$^*$, O.V.Ogievetsky$^{\dagger}$ and P.N.Pyatov$^{*}$}
\\[2em]
{\small $^*$ \it Bogoliubov Laboratory of Theoretical Physics, JINR,
141980  Dubna,} \\ {\small \it Moscow region, Russia}
\vskip .4cm
{\small $^{\dagger}$ \it Center of Theoretical Physics, Luminy,
13288 Marseille, France} \\
{\small and} \\
{\small \it P. N. Lebedev Physical Institute, Theoretical Department,
Leninsky pr. 53, 117924 Moscow, Russia} \\
\end{center}

\vskip 1cm
\centerline{{\bf Abstract}}
\vskip .3cm
\noindent
In the framework of the Drinfeld theory of twists in
Hopf algebras we construct quantum matrix algebras which generalize
the Reflection Equation and the RTT algebras.
Finite-dimensional representations of these algebras
related to the theory of nonultralocal spin chains are
presented. The Cayley-Hamilton-Newton identities
are demonstrated. These identities
allow to define the quantum spectrum for
the quantum matrices. We mention possible applications
of the new quantum matrix algebras to constructions of
noncommutative analogs of Minkowski space and quantum
Poincar\'e algebras.

\vskip 1cm

\section{Twisted Hopf Algebras and related
Quantum Matrix Algebras}

Consider a quasitriangular Hopf algebra
${\cal A}$ $(\Delta, \, \epsilon, \, S, \, {\cal R})$
with a universal $R$-matrix ${\cal R} \in {\cal A} \otimes {\cal A}$
\be
\lb{1}
\Delta' (a) = {\cal R} \, \Delta(a) \, {\cal R}^{-1} \; .
\ee
Let ${\cal F}$
be an invertible element of ${\cal A} \otimes {\cal A}$
$$
{\cal F} = \sum_i \alpha_i \otimes \beta_i \; , \;\;\;
(\epsilon \otimes id) {\cal F} = 1 =
(id \otimes \epsilon ) {\cal F} \; , \;\;\;
{\cal F}^{-1} = \sum_i \, \gamma_i \otimes \delta_i
$$
which defines a twist \cite{D4}
of the algebra ${\cal A}$ to a new quasitriangular Hopf algebra
${\cal A}^{(F)}$
$(\Delta^{(F)}, \, \epsilon^{(F)}, \, S^{(F)}, \, {\cal R}^{(F)})$.
Here $\epsilon^{(F)} = \epsilon$,
a new coproduct $\Delta^{(F)}$ and a new antipode $S^{(F)}$ are
\be
\lb{2}
\Delta^{(F)}(a) = {\cal F} \, \Delta(a) \, {\cal F}^{-1} \; ,
\;\;\; S^{(F)}(a) = u \, S(a) \, u^{-1} \; , \;\;\;
\forall a \in {\cal A}
\ee
and a new universal $R$-matrix is
\be
\lb{2a}
{\cal R}^{(F)} = {\cal F}^{21} \, {\cal R} \, {\cal F}^{-1} \; .
\ee
One can derive the following expressions for the element $u$ in (\ref{2})
$$
u = \sum_i \, \alpha_i \, S(\beta_i) \; , \;\;\;
u^{-1} = \sum_i \, S(\gamma_i) \, \delta_i \; , \;\;\;
S(\alpha_i) \, u^{-1} \, \beta_i = 1 \; .
$$
The coassociativity condition for $\Delta^{(F)}$ is implied by
\be
\lb{3}
{\cal F}^{12} \, (\Delta \otimes id) {\cal F} =
{\cal F}^{23} \, (id \otimes \Delta) {\cal F} \; .
\ee
Impose additional relations on ${\cal F}$
\be
\lb{4}
(\Delta \otimes id) {\cal F} =
{\cal F}^{13} \, {\cal F}^{23} \;\; , \;\;\;
(id \otimes \Delta) {\cal F} =
{\cal F}^{13} \, {\cal F}^{12} \; ,
\ee
which, together with (\ref{3}), imply the Yang-Baxter equation for
${\cal F}$.
Using (\ref{1}) one deduces from (\ref{4}) equations
\be
\lb{5}
{\cal R}^{12} \, {\cal F}^{13} \, {\cal F}^{23} =
{\cal F}^{23} \, {\cal F}^{13} \, {\cal R}^{12} \;\; , \;\;\;
{\cal F}^{12} \, {\cal F}^{13} \, {\cal R}^{23} =
{\cal R}^{23} \, {\cal F}^{13} \, {\cal F}^{12} \; .
\ee
With the Yang-Baxter relations
\be
\lb{6}
{\cal R}^{12} \, {\cal R}^{13} \, {\cal R}^{23} =
{\cal R}^{23} \, {\cal R}^{13} \, {\cal R}^{12} \;\; , \;\;\;
{\cal F}^{12} \, {\cal F}^{13} \, {\cal F}^{23} =
{\cal F}^{23} \, {\cal F}^{13} \, {\cal F}^{12} \; .
\ee
eqs. (\ref{5}) define the twist which is proposed in
\cite{Resh1} (we do not assume that
${\cal F}^{21}{\cal F} = 1 \otimes 1$ as in \cite{Resh1}).
We call the pair of $({\cal R}, \, {\cal F})$
{\it compatible} if eqs (\ref{5}) and (\ref{6}) are satisfied.
Many explicit examples of {\it compatible}
pairs are known (see e.g. \cite{KLM} and references therein).
Using relations (\ref{5}), (\ref{6})
one proves an identity
\be
\lb{6a}
{\cal R}^{12} \, ({\cal R}^{13} \, {\cal F}^{31}) \,
{\cal F}^{21} \, ({\cal R}^{23} \, {\cal F}^{32}) \,
({\cal F}^{21})^{-1} =
({\cal R}^{23} \, {\cal F}^{32}) \, {\cal F}^{12} \,
({\cal R}^{13} \, {\cal F}^{31}) \, ({\cal F}^{12})^{-1} \,
{\cal R}^{12} \; .
\ee
Consider
the dual Hopf algebra ${\cal A}^{*}$
(the algebra of linear functionals on ${\cal A}$)
with generators
$\{ T^{i}_{j} \}$ $(i,j = 1, \dots ,n)$ and the coproduct
$\Delta(T^{i}_{j})=
T^{i}_{k} \otimes T^{k}_{j}$.
The matrix representations
of the algebra ${\cal A}$
in $n$-dimensional vector space $V$
is given by $T^i_j (a) = < a, T^i_j >$
$\forall a \in {\cal A}$ where $<,>$ is the pairing between
${\cal A}$ and ${\cal A}^{*}$.
Introduce the
quantum matrices $L^{\pm} \in {\cal A}$ and $K^{\pm} \in {\cal A}$
\be
\lb{7}
\begin{array}{c}
(L^{+})^{i}_{j} :=  < {\cal R} , id \otimes T^i_j > \; , \;\;
S((L^{-})^{i}_{j}) :=  < {\cal R} , T^i_j \otimes id> \; , \;\; \\[1em]
(K^{+})^{i}_{j} :=  < {\cal F} , id \otimes T^i_j > \; , \;\;
S((K^{-})^{i}_{j}) :=  < {\cal F} , T^i_j \otimes id> \; ,
\end{array}
\ee
and numerical {\it compatible} matrices $R$ and $F$
\be
\lb{8}
\begin{array}{c}
R_{12} = < {\cal R} , T_1 \otimes T_2> =
< L^{+}_2 , T_1> = < S(L^{-}_1) , T_2>
\;\; , \\[1em]
F_{12} = < {\cal F} , T_1 \otimes T_2> =
<K^{+}_2 , T_1> =
<S(K^{-}_1) , T_2> \; .
\end{array}
\ee
Matrices $R_{12}$ and $F_{12}$ act in $V \otimes V$.
Further, we denote  the matrices $F$ and $R$ acting in
$V_k \otimes V_{k+1}$ by $F_k$ and $R_k$
(here the subscripts $k$ and $k+1$ enumerate different copies of the space $V$).

Our aim is to define a quantum matrix
algebra which generalizes the $RTT$ \cite{FRT}
and the Reflection Equation algebras \cite{KSa}.
To this end, define a quantum matrix $M^i_j \in {\cal A}$
\be
\lb{13}
M^i_j :=
< {\cal R}^{21} \, {\cal F}^{12} , id \otimes T^i_j > =
S((L^{-})^i_k) \, (K^{+})^k_j  \; .
\ee
and consider the following pairing
\be
\lb{14}
\begin{array}{c}
< {\cal R}^{21} \, {\cal F}^{12} , id \otimes T_1 \, T_2 > =
< {\cal R}^{21} \, {\cal R}^{31} \, {\cal F}^{13} \, {\cal F}^{12}
, id \otimes T_1 \otimes T_2 > = \\[1em]
S(L^{-}_1) \, S(L^{-}_2) \, (K^{+}_2) \, K^{+}_1 =
S(L^{-}_1) \, S(L^{-}_2) \, \hat{F} \, (K^{+}_2) \, K^{+}_1 \,
\hat{F}^{-1} = \\[1em]
S(L^{-}_1) \, K^{+}_1 \, \hat{F} \,
S(L^{-}_1) \, K^{+}_1 \, \hat{F}^{-1} =
M_1 \hat{F} M_1 \hat{F}^{-1} \; .
\end{array}
\ee
Here $\hat{F} = P_{12} \, F_{12}$.
The pairing of the identity (\ref{6a}) with
$( T_1 \otimes T_2 \otimes id )$
gives (with the help of (\ref{14}))
the commutation relations for
the elements of the matrix $M$:
\be
\lb{15}
\R \, M_1 \, \hat{F} \, M_1 \, \hat{F}^{-1} =
M_1 \, \hat{F} \, M_1 \, \hat{F}^{-1} \, \R \; ,
\ee
where $\R:=P_{12}R_{12}$.
The algebras ${\cal M}(\R,\F)$ with generators $M^i_j$ and defining
relations
(\ref{15}) unify the $RTT$ ($\F=P$) and the Reflection Equation
($\F = \R$) algebras. The algebras of that kind
have been considered in \cite{Hlav} (see also \cite{KSa},
\cite{KSk}).

The comultiplication acts on $M^i_j$ as follows
\be
\lb{15a}
\Delta(M^i_j) = M^l_k \otimes S((L^{-})^i_l) \, (K^+)^k_j \; ,
\ee
and one can
define for the algebra ${\cal M}(\R,\F)$ (\ref{15}) the homomorphic
mapping
\be
\lb{braid}
\overline{\Delta}(M^i_j) = M^i_k \, \widetilde{M}^k_j
\; , \;\;\;
\hat{F} \, M_1 \, \hat{F}^{-1} \, \widetilde{M}_1 =
\widetilde{M}_1 \, \hat{F} \, M_1 \, \hat{F}^{-1} \; ,
\ee
which is an analog of the braided coproduct \cite{Ma}
(elements $\widetilde{M}^i_j$
generate another copy of the algebra ${\cal M}(\R,\F)$ (\ref{15})).

The extension of (\ref{14}) on arbitrary products of $T$'s
is straightforward
\be
\lb{16}
< {\cal R}^{21} \, {\cal F}^{12} , id \otimes T_1 \, T_2 \dots T_n > =
M_{\overline{1}} \, M_{\overline{2}} \dots M_{\overline{n}} \; ,
\ee
where
\be
\lb{16a}
M_{\overline{1}} = M_1 \; , \;\;\;
M_{\overline{k+1}} = (\hat{F}_k M_{\overline{k}} \hat{F}^{-1}_k) \; .
\ee

Comparing relations (\ref{14}) and (\ref{15}) one can
conclude that ${\cal A}^*$ is the $RTT$ algebra
with defining relations
\be
\lb{RTT}
\R \, T_1 \, T_2 = T_1 \, T_2 \, \R \; ,
\ee
if we require that the
linear mapping $\phi : \;\; {\cal A}^* \rightarrow {\cal A}$
\be
\lb{map}
\phi(\alpha) :=
< {\cal R}^{21} \, {\cal F}^{12} , id \otimes \alpha >
\; \in \; {\cal A} \;\; , \;\;\; \alpha \in {\cal A}^* \; ,
\ee
is nondegenerate: $\phi(\alpha)=0 \Rightarrow \alpha =0$.

\section{Generalized Cayley-Hamilton-Newton \\
Identities}

Consider the case when $\R_{12}$ (\ref{8})
is a Hecke type matrix
\be
\lb{Heck}
\R^2 = (q-q^{-1}) \, \R + I \; ,
\ee
where $I$ is the identity matrix and $q$ is a numerical
parameter which is not equal to a root of unity.
In this case we have proved \cite{IOP2} that the $q$-matrix $T$
(\ref{RTT}) satisfies
the generalized Cayley-Hamilton-Newton
identities
\be
\lb{chnT}
i_q T^{\overline{\wedge i}} =
  \sum_{k=0}^{i-1}(-1)^{i-k+1} T^{\overline{i-k}}\ \sigma_k(T) \; ,
\;\;\;
i_q T^{\overline{{\cal S} i}} =
  \sum_{k=0}^{i-1} \, T^{\overline{i-k}}\ \tau_k(T) \; ,
\ee
where $k_q = (q^k - q^{-k})/(q-q^{-1})$,
\be
\lb{def1}
T^{\overline{\wedge i}} :=
\Rtr_{(2\dots k)} \left( A^{(k)} \, T_1 \dots T_k \right)
\; , \;\;\;
T^{\overline{{\cal S} i}} :=
\Rtr_{(2\dots k)} \left( S^{(k)} \, T_1 \dots T_k \right)
\; ,
\ee
are versions of $k$-wedge and $k$-symmetric powers of
$q$-matrix $T$ \cite{IOP2},
\be
\lb{def1a}
\sigma_k(T) := q^k \,
\Rtr_{(1\dots k)} \left( A^{(k)} \, T_1 \dots T_k \right)
\; , \;\;\;
\tau_k(T) := q^{-k} \,
\Rtr_{(1\dots k)} \left( S^{(k)} \, T_1 \dots T_k \right)
\; ,
\ee
($\sigma_0(T) = \tau_0(T)= 1$)
are elementary and complete symmetric functions of the
spectrum of the $q$-matrix $T$ \cite{IOP2},
\be
\lb{def2}
T^{\overline{k}} :=
\Rtr_{(2\dots k)}
\left( \R_1 \dots \R_{k-1} \, T_1 \dots T_k \right) \; , \;\;\;
T^{\overline{1}} := T \; ,
\ee
are generalized $k$-th power of the $q$-matrix $T$
\cite{IOP1}. Projectors $A^{(k)}$ and  $S^{(k)}$ are
$q$-antisymmetrizers
and $q$-symmetrizers respectively
\ba
\lb{antis}
A^{(1)} := \id\  ,&\quad&
A^{(k)} := {1\over k_q}\,
A^{(k{-}1)}\left(q^{k-1}-(k{-}1)_q\R_{k{-}1}\right)A^{(k{-}1)}\ ,
\\
\lb{simm}
S^{(1)} := \id\  ,&\quad&
S^{(k)} := {1\over k_q}\,
S^{(k{-}1)}\left(q^{1-k}+(k{-}1)_q\R_{k{-}1}\right)S^{(k{-}1)}\ .
\ea
The symbol $\Rtr$ in (\ref{def1})-(\ref{def2}) denotes
the quantum trace $\Rtr(Z) := \tr(D \, Z)$ where a numerical
matrix $D$ is defined by the $R$-matrix
\be
\lb{rD}
D = \tr_{(2)} \Psi_{12} \; , \;\;\;
\tr_{(2)} ( \Psi_{12} \, \R_2) = P_{13} =
\tr_{(2)} (\R_1 \, \Psi_{23} ) \; .
\ee
The following properties of matrix $D$ will be important
\be
\lb{DR}
\Rtr_{2} \R_1 = I_1 \; , \;\;\;
\Rtr_{2} (\R^{\pm 1}_1 \, Z_1 \, \R_1^{\mp 1}) =
I_1 \, \Rtr (Z) \; \Rightarrow \;
\R_1 \, D_1 \, D_2 = D_1 \, D_2 \, \R_1 \; ,
\ee
\be
\lb{DF}
\Rtr_{2} (\F^{\pm 1}_1 \, Z_1 \, \F_1^{\mp 1}) =
I_1 \, \Rtr (Z) \; \Rightarrow \;
\F_1 \, D_1 \, D_2 = D_1 \, D_2 \, \F_1
\; ,
\ee
where $Z$ is an arbitrary quantum $(n \times n)$ matrix.
One can verify these properties by using matrix versions
of the relations (\ref{5}), (\ref{6}). Note
that the identities (\ref{chnT}) and the definitions (\ref{def2})
have been written in \cite{IOP2} in a different form
with the usual traces instead of the quantum ones.
It is clear, however, that these
two forms of Cayley-Hamilton-Newton identities are
equivalent since they are
related by a simple automorphism of the $RTT$ algebra
$T \rightarrow D \, T$.

Now we show that for the case of algebras ${\cal M}(\R,\F)$ (\ref{15})
(unifying the reflection equation and the $RTT$ algebras)
one can also prove the identities which
are analogous to the
generalized Cayley-Hamilton-Newton identities (\ref{chnT}).

\vspace{0.5cm}
\noindent
{\bf Theorem} {\it
The generalized Cayley-Hamilton-Newton
identities for the algebra ${\cal M}(\R,\F)$
with defining relations (\ref{15}) have the form
\be
\lb{chnM}
i_q \, M^{\wedge i} =
  \sum_{k=0}^{i-1} \, (-1)^{i-k+1} \, M^{\overline{i-k}}\ \sigma_k(M)
\; , \;\;\;
i_q \, M^{{\cal S} i} =
  \sum_{k=0}^{i-1} \, M^{\overline{i-k}}\ \tau_k(M)\ ,
\ee
where
\be
\lb{def11}
\begin{array}{c}
M^{\overline{k}} :=
\Rtr_{(2\dots k)}
\left( \R_1 \dots \R_{k-1} \, M_{\overline{1}} \dots M_{\overline{k}}
\right) \; , \\[1em]
M^{\overline{\wedge i}} :=
\Rtr_{(2\dots k)} \left( A^{(k)} \, M_{\overline{1}} \dots
M_{\overline{k}}
 \right) \; , \\[1em]
M^{\overline{{\cal S} i}} :=
\Rtr_{(2\dots k)} \left( S^{(k)} \, M_{\overline{1}} \dots
M_{\overline{k}}
\right) \; ,
\end{array}
\ee
are quantum versions of the $k$-th, $k$-wedge and $k$-symmetric powers
of
the $q$-matrix $M$,
\be
\lb{def11a}
\begin{array}{c}
\sigma_k(M) := q^k \,
\Rtr_{(1\dots k)} \left( A^{(k)} \, M_{\overline{1}} \dots
M_{\overline{k}}
\right)
\; , \\
\tau_k(M) := q^{-k} \,
\Rtr_{(1\dots k)} \left( S^{(k)} \, M_{\overline{1}} \dots
M_{\overline{k}}
 \right)
\; ,
\end{array}
\ee
are elementary and complete symmetric functions of
the spectrum of the $q$-matrix $M$.
The elements $\sigma_k(M)$, $\tau_k(M)$ generate a commutative subalgebra in
${\cal M}(\R,\F)$.
}

\vspace{0.5cm}
\noindent
{\bf Proof.}
The identities (\ref{chnM}) can be obtained directly
by the map (\ref{map}) from the $RTT$ identities (\ref{chnT}).
To do this, consider the map
$$
< {\cal R}^{21} \, {\cal F}^{12} , id \otimes
T^{\overline{i-k}}\ e_k(T) > =
< {\cal R}^{21} \, {\cal R}^{31} \, {\cal F}^{13} \, {\cal F}^{12} \, ,
id \otimes T^{\overline{i-k}} \otimes e_k(T) > =
$$
$$
= < {\cal R}^{21} \, , id \otimes (T^{\overline{i-k}})_{(1)} > \,
< {\cal R}^{31} \, {\cal F}^{13} \, , id \otimes 1 \otimes e_k(T) > \,
<  {\cal F}^{12} \, , id \otimes (T^{\overline{i-k}})_{(2)} > =
$$
\be
\lb{17}
= f_{(1)}(S(L^-)) \, e_k(M)  \, f_{(2)}(K^+) \; ,
\ee
where $e_k(M) = \sigma_k(M), \, \tau_k(M)$,
$$
\Delta  (T^{\overline{i-k}}) := (T^{\overline{i-k}})_{(1)}
 \otimes (T^{\overline{i-k}})_{(2)}  \; .
$$
and $f_{(1,2)}(.)$ are some functions
explicit forms of which are not important.
We need only the relation $[ e_k(M) , f_{(2)}(K^{+}) ] =0$
which is deduced from the relation
$[e_k(M), \,  K^{+} ] = 0$.
Indeed,
$$
\begin{array}{c}
\Rtr_{(1,\dots ,k)}
(A^{(1,k)} M_{\overline{1}}\dots M_{\overline{k}}) \, K^+_0 = \\[1em]
K^+_0 \,
\Rtr_{(1,\dots ,k)} (\F_{0 \rightarrow k-1} \,
A^{(0,k-1)} M_{\overline{0}}\dots M_{\overline{k-1}}
\F^{-1}_{0 \rightarrow k-1})
\end{array}
$$
\be
\lb{18}
= K^+_0 \,
\Rtr_{(0,\dots ,k-1)} \,
(A^{(0,k-1)} M_{\overline{0}}\dots M_{\overline{k-1}} ) =
K^{+}_0 e_k(M) \; .
\ee
Here we have used $M_1 K^{+}_0 = K^{+}_0 \, \F_{0} \, M_0 \,
\F^{-1}_{0}$.

Using (\ref{18}) one can rewrite (\ref{17}) in the form
\be
\lb{19}
< {\cal R}^{21} \, {\cal F}^{12} , id \otimes
T^{\overline{i-k}}\ e_k(T) > =
< {\cal R}^{21} \, {\cal F}^{12} , id \otimes
T^{\overline{i-k}} > e_k(M) =
M^{\overline{i-k}} \, e_k(M) \; .
\ee
Now it is evident that identities
(\ref{chnT}) and (\ref{chnM}) are related via the map (\ref{16}).

The commutativity of the elements $\sigma_k(M)$ and $\tau_k(M)$b is proved
by the method presented in \cite{IOP3}.

\vspace{0.5cm}
\section{Automorphisms for the algebras
${\cal M}(\R,\F)$.}

In the paper \cite{IOP3} we have obtained the
quantum
Cayley-Hamilton-Newton identities for the algebra
of $q$-matrices $M^i_j$ with different
defining relations (cf. with (\ref{15}))
\be
\lb{21}
\R \, M_1 \, \hat{F} \, M_1 \, \hat{F}^{-1} =
M_1 \, \hat{F} \, M_1 \, \hat{F}^{-1} \, \R^{\F \F} \; .
\ee
where $\R$ is a Hecke type $R$- matrix,
the pair $(\R, \, \F)$ is {\it compatible} and \cite{IOP3}
\be
\lb{rff}
\R^{\F \F} := \F^2 \, \R \, \F^{-2} =
D_1' \, D_2' \, \R (D_1' \, D_2')^{-1} \; .
\ee
Here a numerical matrix $D'$ is an analog of the matrix $D$
(\ref{rD}):
\be
\lb{fD}
D' = \tr_{(2)} \Phi_{12} \; , \;\;\;
\tr_{(2)} ( \Phi_{12} \, \F_2) = P_{13} =
\tr_{(2)} (\F_1 \, \Phi_{23} ) \; ,
\ee
but is related to the matrix $\F$.

The remarkable fact is that these two algebras
${\cal M}(\R,\F)$ (\ref{15}) and (\ref{21})
are connected by a simple transformation.
To manifest this connection we need to introduce
an algebra which unifies the both algebras
(\ref{15}) and (\ref{21}). This algebra
is defined by relations
\be
\lb{22}
\R \, M_1 \, \hat{A} \, M_1 \, (\F^{-1})^X =
M_1 \, \hat{A} \, M_1 \, (\F^{-1})^{X} \, \R^{X} \; .
\ee
where $(\F^{-1})^X = X_1 X_2 \F^{-1} (X_1 X_2)^{-1}$,
$\R^X = X_1 X_2 \R (X_1 X_2)^{-1}$, $X^i_j$ is
an invertible numerical $(n \times n)$ matrix,
the pair $(\R, \, \F)$ is {\it compatible} and
matrix $\hat{A}$ is defined by the matrix $\Phi$,
the skew inverse of the matrix $\F$
(\ref{fD})
\be
\lb{23}
\hat{A}_{12} := X_2 \, \Phi^{-1}_{21} \, X^{-1}_1 \; .
\ee
Thus, the algebra (\ref{22}) is defined by the {\it compatible} pair
$(\R, \, \F)$ and the invertible matrix $X$. For a
Hecke-type $\R$ and an arbitrary $X$ one can prove,
using the methods of the paper \cite{IOP3}, the
Cayley-Hamilton-Newton identities similar to (\ref{chnM})
for the algebra (\ref{22}).

For $X=D$ we find that $(\F^{-1})^X = \F^{-1}$, $\R^X = \R$,
$\hat{A} = \F$ (these equalities can be obtained from
(\ref{DF}) and (\ref{23})) and therefore the algebra (\ref{22})
coincides with the algebra ${\cal M}(\R,\F)$ (\ref{15}).

On another hand, for $X=D'$ we have $(\F^{-1})^X = \F^{-1}$,
$\R^X = \R^{\F \F}$ (\ref{rff}),
$\hat{A} = \F$ and, thus, the algebra (\ref{22})
converts to the algebra (\ref{21}).

Now we note that relations (\ref{22}) are invariant under
the transformations
\be
\lb{24}
\begin{array}{c}
\R \longrightarrow U_1 \, U_2 \, \R \, (U_1 U_2)^{-1} \; , \;\;\;
\F \longrightarrow U_1 \, U_2 \, \F \, (U_1 U_2)^{-1} \; , \\[2em]
\hat{A} \longrightarrow Y_1 \, U_2 \,
\hat{A} \, (U_1 Y_2)^{-1} \; ,  \;\;\;
M \longrightarrow U \, M \, Y^{-1} \; , \;\;\;
X \longrightarrow Y \, X \, U^{-1} \; ,
\end{array}
\ee
where $U$ and $Y$ are invertible numerical $(n \times n)$
matrices. One can always find such matrices $U,Y$ that
$X = D \rightarrow X = D'$ and
therefore relate the algebras
of the types (\ref{15}) and (\ref{21}) by
transformations (\ref{24}). So, using (\ref{24})
one can obtain the
Cayley-Hamilton-Newton identities for one algebra from another
and vise versa.

\vspace{0.5cm}
\section{Discussion}

\noindent
{\bf 1.} The generators $M^i_j$ of the matrix
algebra ${\cal M}(\R,\F)$ (\ref{15})
have a natural finite-dimensional matrix representation
$$
(M^i_j)^n_m = (\R \, \F)^{ni}_{mj} \; ,
$$
which is dictated by the definition (\ref{13}).
This representation
and the braided coproduct $\overline{\Delta}(M^i_j)$
(\ref{braid}) give a
possibility to define the integrable nonultralocal
spin chains following
the approach of the paper \cite{Hlav}.

\noindent
{\bf 2.} The explicit form (\ref{15a}) for the
comultiplication $\Delta(M^i_j)$
leads to the conclusion that
the algebra (\ref{22})
(and, hence, its particular formulations (\ref{15}) and (\ref{21})) are covariant
under the quantum transformations
\be
\lb{25}
M^i_j \rightarrow (T \, M \, \widetilde{T}^{-1})^i_j \equiv
M^k_l \otimes T^i_k \, S(\widetilde{T}^l_j)
\ee
where elements $\{ T^i_j , \, \widetilde{T}^k_l \}$
generate the algebra
\be
\lb{25a}
\begin{array}{c}
\R \, T_1 \, T_2 =  T_1 \, T_2 \, \R \; , \;\;\;
\hat{A} \, T_1 \, \widetilde{T}_2 =
\widetilde{T}_1 \, T_2 \, \hat{A} \; , \\[1em]
\R^X \, \widetilde{T}_1 \, \widetilde{T}_2 =
\widetilde{T}_1 \, \widetilde{T}_2 \, \R^X \; , \;\;\;
\F^X \, \widetilde{T}_1 \, \widetilde{T}_2 =
\widetilde{T}_1 \, \widetilde{T}_2 \, \F^X \;  .
\end{array}
\ee

\noindent
{\bf 3.} There exists a Heisenberg double
of the ${\cal M}(\R,\F)$ (\ref{22}) and
the RTT algebras:
\be
\lb{26}
\begin{array}{c}
\R \, M_1 \, \hat{A} \, M_1 \, \hat{B} =
M_1 \, \hat{A} \, M_1 \, \hat{B} \, \R \; , \\[1em]
T_1 \, M_2 = \hat{A} \, M_1 \, \hat{B} \, T_1 \; , \\[1em]
\R \, T_1 \, T_2 = T_1 \, T_2 \, \R \; .
\end{array}
\ee
where we rewrite the relations (\ref{22})
using $\hat{B} = X_1 X_2 \F^{-1}$.
For this Heisenberg double one can define the analog
of Alekseev-Faddeev discrete evolution (see \cite{AF}) $M \rightarrow M$,
$T \rightarrow T^{(k)} = M^k \, T$.

\noindent
{\bf 4.} One can define the de Rham differential complex
over the algebra ${\cal M}(\R,\F)$ (\ref{22}) (for matrix $\R$
of Hecke type)
\be
\lb{27}
\begin{array}{c}
\R \, M_1 \, \hat{A} \, M_1 \, \hat{B} =
M_1 \, \hat{A} \, M_1 \, \hat{B} \, \R \; , \\[1em]
\R \, (d M_1) \, \hat{A} \, M_1 \, \hat{B} =
M_1 \, \hat{A} \, (d M_1) \, \hat{B} \, \R^{-1} \; , \\[1em]
\R \, (d M_1) \, \hat{A} \, (d M_1) \, \hat{B} =
- (d M_1) \, \hat{A} \, (d M_1) \, \hat{B} \, \R^{-1} \; .
\end{array}
\ee

\noindent
{\bf 5.}
The construction of the Heisenberg double
(\ref{26}) and the covariance (\ref{25}) give a
possibility to apply the algebra
${\cal M}(\R,\F)$ (\ref{21}) in the $(2 \times 2)$
case to the theory of a quantum Minkowski space and a quantum
Poincar\'e algebra \cite{OSWZ}, \cite{AKR}.
We note here that this subject in the framework
of twist theory has been already discussed  \cite{KM}
in different context.

\noindent
{\bf 6.} The algebra ${\cal M}(\R , \F)$ in the form (\ref{21})
for the dynamical $\R$ and $\F$ matrices has been used
for $R$-matrix quantization of the
elliptic Ruijsenaars-Schneider model \cite{ACF}.

\vskip 1cm

\noindent
{\bf Acknowledgements:} We are indebted to
L.Hlavaty, P.~Kulish and S.Pakuliak
for valuable discussions. This work was supported
in parts by CNRS grant PICS No. 608 and RFBR grant No.
98-01-2033. The work of AI and PP was also supported by
the RFBR grant No. 97-01-01041.

\end{document}